\documentclass[preprint,12pt]{elsarticle}

\usepackage{amssymb}

\usepackage{bm}
\usepackage{amsmath}
\usepackage{amsfonts}

\def\R{\mathbb{R}}

\def\Z{\mathbb{Z}}
\def\N{\mathbb{N}}

\newcommand\etal{\mbox{\textit{et al.}}}

\journal{Applied Numerical Mathematics}

\begin{document}

\begin{frontmatter}



\title{Analysis and discretization of the volume penalized Laplace operator with Neumann boundary conditions}


\author[label1]{Dmitry Kolomenskiy}

\address[label1]{Department of Mathematics and Statistics, McGill University, CRM, Montr\'eal, Canada}

\author[label2]{Romain Nguyen van yen}

\address[label2]{FB Mathematik und Informatik, Freie Universit\"at Berlin, Germany}

\author[label3]{Kai Schneider}

\address[label3]{M2P2-CNRS, Aix-Marseille Universit\'e, 38 rue Joliot-Curie, 13451 Marseille Cedex 20, France}

\begin{abstract}
We study the properties of an approximation of the Laplace operator with Neumann
boundary conditions using volume penalization.
For the one-dimensional Poisson equation we compute explicitly the exact solution of the penalized equation
and quantify the penalization error.
Numerical simulations using finite differences 
allow then to assess the discretisation and penalization errors.
The eigenvalue problem of the penalized Laplace operator with Neumann boundary conditions is also studied.
As examples in two space dimensions, we consider a Poisson equation with Neumann boundary conditions in rectangular and circular domains.
\end{abstract}


\end{frontmatter}


\section{Introduction}

Solving partial differential equations (PDEs) in complex domains is unavoidable in real world
applications. Different numerical methods have been developed so far, for example body fitted computational grids
or coordinate transforms \cite{FePe96}. 
Immersed boundary methods are still of growing interest due to their high flexibility
and their ease of implementation into existing codes. 
The underlying idea of these methods is to embed the complex geometry into a simple geometry ({\it e.g.} a rectangle)
for which efficient solvers are available. The boundary conditions are then imposed by adding supplementary terms
to the governing  equations. Different penalization approaches are on the market, for example,
surface and volume penalization techniques, immersed boundary methods using direct forcing
and Lagrangian multipliers.
For reviews on immersed boundary techniques, we refer to~\cite{Pesk02,MiIc05}. 

In the current work, we focus on the volume penalization method developed by Angot \etal~\cite{ABF99} for
imposing Dirichlet boundary conditions in viscous fluid flow. Physically, the boundary conditions correspond to no-slip
conditions on the wall, {\it i.e.}, both the normal and the tangential velocity do vanish at the fixed wall.
This penalization approach is physically motivated as walls or solid obstacles are modeled as porous media
whose permeability tends to zero.
Mathematically, it has also been justified. In \cite{ABF99, CaFa03} it was shown that the solution of the penalized Navier--Stokes equations
converges towards the solution of the Navier--Stokes equations with no-slip boundary conditions, while the error depends on the penalization parameter.
Various applications of the volume penalization method to impose Dirichlet boundary conditions
can be found in the literature. Briefly summarizing, we can mention 
computations of confined hydrodynamic and magnetohydrodynamic turbulence, which can be found in \cite{ScFa05} and \cite{SNB11,MLBS12}, respectively.
Fluid-structure interaction simulations have been carried out for moving obstacles \cite{KoSc09}
and for flexible beams \cite{KoES13}. Applications to the aerodynamics of insect flight in two
and three space dimensions can be found in \cite{KMFS11}.

Most of the developed penalization techniques deal with Dirichlet boundary conditions, and only few allow to impose Neumann conditions. 
Neumann boundary conditions in partial differential equations are encountered in many applications,
for example when solving the Poisson equation for pressure in incompressible flows, to model adiabatic walls in heat transfer, 
or to impose no-flux conditions for passive or reactive scalars at walls.
In \cite{BoLe05} a review on the pure Neumann problem using finite elements is given
and different techniques for solving the algebraic system are discussed.
An extension of the volume penalization method \cite{ABF99} to impose Neumann or Robin boundary conditions has been presented
in \cite{RaB07} and applied in the context of finite element or finite volumes~\cite{RAB07}.
In \cite{KKAS12} we extended this method for pseudo-spectral discretizations
and applied it to scalar mixing in incompressible flow for fixed and also for moving geometries
imposing no-slip conditions for the velocity and no-flux conditions for the passive scalar field.

The fields of possible applications of the volume penalization method for imposing Neumann conditions in complex geometries are multifarious and large.
For example, confined magnetohydrodynamic flow configurations can be studied imposing finite values of the current density at the wall, or convection problems which necessitate imposing a given heat flux at the boundary.

Motivated by the work of \cite{MiGo03}, where
the properties of Fourier approximations of elliptic problems with discontinuous coefficients have been studied,
we analyzed mathematically the penalized Laplace and Stokes operators
with Dirichlet boundary conditions in \cite{NKS12} and verified the predicted convergence numerically.
%
The aim of the present work is to generalize the approach developed in \cite{NKS12}
and to analyze the penalized Laplace operator with Neumann boundary conditions.
For a one-dimensional Poisson equation, we explicitely compute the penalization error by solving the penalized equation 
analytically. Discretizing the penalized equation using finite difference methods,
we study the influences of both the numerical resolution and the value of the penalization parameter.

The outline of the paper is the following:
First we consider the penalized Poisson equation in one space dimension with Neumann boundary conditions both analytically and numerically.
Then, in section~3 we study the eigenvalue problem of the penalized Laplace operator.
Section~4 presents applications of the penalization method to solve the Poisson equation in two dimensions in 
a rectangular and a circular domain.
Finally, some conclusions are drawn and some perspectives are given in section~5.

\section{Poisson equation with Neumann boundary conditions and penalization}

\subsection{Problem setting}

We consider the one-dimensional Poisson equation
\begin{equation}   
- w'' \, = \, f \quad \mbox{\rm for} \quad x \in (0, \pi)
\label{poisson1d}
\end{equation}
completed with homogeneous Neumann boundary conditions, $w'(x=0) = w'(x=\pi) = 0$ and for $f(x) = m^2 \cos m x$, $m \in \Z$.
The exact solution $w \in H^2(0, \pi)$ is given by $w(x) = \cos mx +C$, where $C \in R$ is an arbitrary constant, 
as the solution is not unique.
Integrating eq.~(\ref{poisson1d}) over $(0, \pi)$ yields the compatibility condition $\int_{0}^{\pi} f(x) dx = w'(x=\pi) - w'(x=0) = 0$
which has to be satisfied to guarantee the existence of a solution.

Following \cite{KKAS12}, the penalized Poisson equation
reads
\begin{equation}
- d_x ((1 - \chi) + \eta \chi) d_x v  \, = \, f \quad \mbox{\rm for} \quad x \in (0, 2 \pi)
\label{poisson1d_penalized}
\end{equation}
where $\eta>0$ is the penalization parameter and $\chi$ the mask function defined by
\begin{equation}
\chi(x) \, = \,   \left \{
    \begin{array} {ll}
        0 \quad \quad \quad \mbox{\rm for} \quad \quad 0 < x < \pi \\
        1/2 \quad \quad \mbox{\rm for} \quad \quad x = 0 \; {\rm or} \; x = \pi \\
        1 \quad \quad \quad \; \mbox{elsewhere}\\
    \end{array}
  \right.
\label{eq:mask1d}
\end{equation}
The domain $\Omega_f = ] 0, \pi [ $, also called fluid domain, is imbedded into the larger domain $\Omega = ] 0, 2 \pi [$ imposing now periodic boundary conditions at the boundary.
Thus we have $\Omega = \Omega_f \cup \Omega_s$, where $\Omega_s$ is the penalization domain, also called solid domain.

\subsection{Analytic solution of the one-dimensional penalized equation}

The penalized Poisson equation (\ref{poisson1d_penalized}) can be solved analytically
in each sub-domain, i.e.,
\begin{eqnarray}
- v'' \, &=& \, f \quad \mbox{\rm for} \quad x \in ]0, \pi [ \\
- \eta v''  \, &=& \, 0 \quad \mbox{\rm for} \quad x \in ]\pi, 2 \pi [
\end{eqnarray}
and accordingly we obtain
\begin{equation}
v(x) \, = \,   \left \{
    \begin{array} {ll}
        \cos m x + A_1 x + A_2 \quad \mbox{\rm for} \quad x \in ]0, \pi [ \\
        B_1 x + B_2            \quad \quad \quad \quad \quad \mbox{\rm for} \quad x \in ]\pi, 2 \pi [
    \end{array}
  \right.
\end{equation}
The coefficients can then be determined by imposing continuity of the solution and of the flux, at $x=0 (= 2\pi)$ and $\pi$,
\begin{eqnarray}
v(\pi^-) \; &=& \; v(\pi^+)         \quad \mbox{\rm and} \quad v(0^+) \; = \; v(2 \pi^-) \\
v'(\pi^-) \; &=& \,  \eta v'(\pi^+) \quad \mbox{\rm and} \quad v'(0^+) \; = \; \eta v'(2 \pi^-)
\end{eqnarray}
This results in
\begin{eqnarray}
A_1 \, &=& \, \frac{1 - (-1)^m}{\pi (1 + 1/\eta)} \quad \mbox{\rm and} \quad  B_1 = \frac{1}{\eta} A_1 \\
A_2 \, &=& \,  \frac{2 \pi}{\eta} \, \frac{1 - (-1)^m}{\pi (1 + 1/\eta)} - 1 + B_2
\end{eqnarray}
Only three of the four coefficients can be determined, $B_2$ corresponds to the additive constant.

%

\begin{figure}[htb]
\centering
\includegraphics{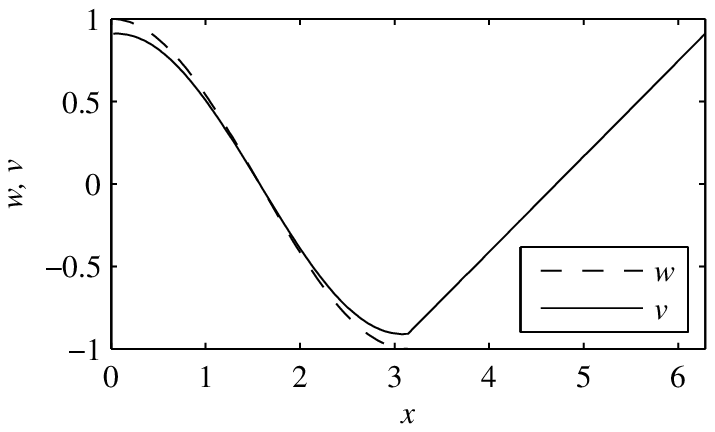}
\includegraphics{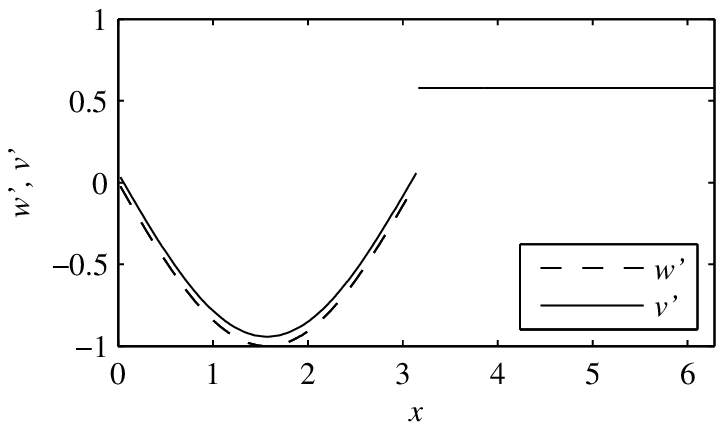}
\includegraphics{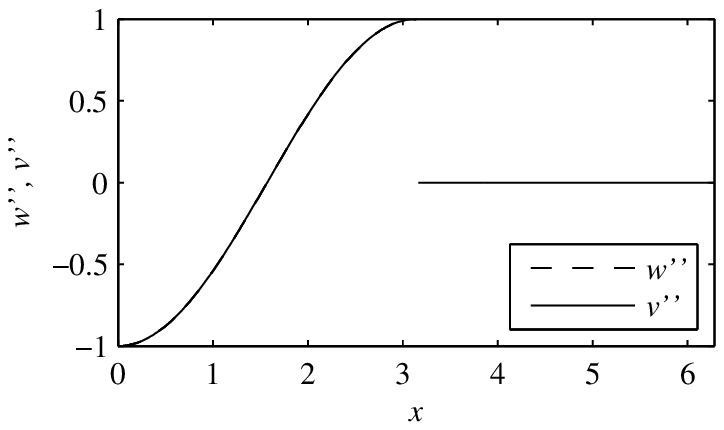}
\caption{Exact solution of the Poisson equation $w(x)$ and exact solution of the penalized equation $v(x)$ using $\eta = 10^{-1}$,
both for $m=1$ (top).
The first (middle) and second (bottom) derivatives are also shown. 
}\label{fig:exact_solution}
\end{figure}

\medskip
Figure~\ref{fig:exact_solution} shows the exact solution, $w(x)$,
and the solution of the penalized problem, $v(x)$ (for $\eta= 10^{-1}$), in the case
$m=1$. Unlike for the penalized heat equation with Neumann boundary conditions \cite{KKAS12}, here there is no boundary
layer in the penalized domain. Note that, if $m$ is even, $v$ and
$w$ coincide exactly. Therefore, in the following let us assume $m$
odd. The coefficients of the penalized solution become (with the
integration constant chosen such as to ensure zero mean value)
\begin{equation}
 A_1 = \frac{2}{\pi} \frac{\eta}{1+\eta}, \quad\quad B_1 = \frac{2}{\pi} \frac{1}{1+\eta}, \quad\quad A_2 = - \frac{\eta}{1+\eta}, \quad\quad B_2 = -
 \frac{3}{1+\eta}.
\label{eq:coeff_modd}
\end{equation}
The difference between the exact solution of the non penalized problem $w$ and $v$ yields the penalization error 
$||w(x) - v(x)||$ which is of order $O({\eta})$ in $\Omega_f$, and which is in this particular case better than the general 
$O(\sqrt{\eta})$ convergence behavior shown in \cite{KKAS12} for the heat equation.

\begin{figure}[htb]
\centering
\includegraphics[width=7.0cm, height=5cm]{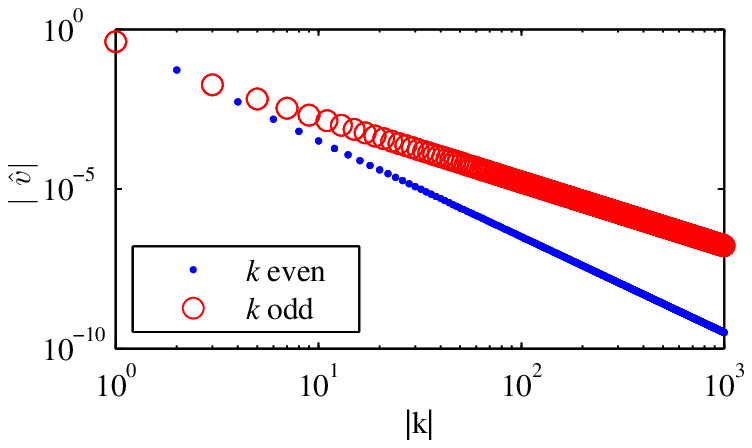}
\caption{Decay of the Fourier coefficients. 
Absolute value of the Fourier coefficients of the exact solution of the penalized equation for $m = 1$ using $\eta = 10^{-1}$. 
The even and odd wavenumbers exhibit different power law behaviors.}\label{fig:fourier_coef}
\end{figure}

It is straightforward to compute the Fourier coefficients of
the solution of the penalized equation $v(x)$:
\begin{equation}
 \hat{v}(k) \, = \, \left\{ \begin{array}{cc}
  \displaystyle \frac{i}{\pi} \frac{m^2}{k(m^2-k^2)} & \space \mathrm{~if~} k~\mathrm{even} \\
  \displaystyle \frac{2}{\pi^2 k^2} \frac{1-\eta}{1+\eta} & \space \mathrm{~if~} k~\mathrm{odd~~and}~~k \ne \pm m \\
  \displaystyle \frac{2}{\pi^2 m^2} \frac{1-\eta}{1+\eta} + \frac{1}{4} & \space \mathrm{~if~} k~\mathrm{odd~~and}~~k = \pm m \\
 \end{array} \right.
\label{eq:fourier_coeff}
\end{equation}
Figure~\ref{fig:fourier_coef} displays the decay of the absolute
value of $\hat{v}$. The leading order is $\sim k^{-2}$ and the
constant pre-factor is finite in the limit $\eta \to 0$. There is
no `intermediate' regime of slow decay at low $k$, because there
is no boundary layer in contrast to the Dirichlet case \cite{NKS12}. 
This rate of decay of $\hat{v}$ suggests
that a Galerkin truncated approximation to $v$ converges as
$N^{-3/2}$. 

\subsection{Discretization error of the second order finite difference scheme}
\label{subsec_defds}

Now we consider the discretization of the penalized equation using centered finite
differences of second order.
Discretizing the equation
\begin{equation}
- d_x (\theta(x)) d_x u  \, = \, f \quad \mbox{\rm for} \quad x \in (0, 2 \pi)
\label{poisson1d_penalized_b}
\end{equation}
where $\theta = (1 - \chi) + \eta \chi$
with periodic boundary conditions
on $N$ grid points $x_i = i/(2 \pi), i = 0, ..., N-1$ yields to the following linear system
\begin{equation}
-D \Theta D = F 
\label{poisson1d_penalized_linsystem}
\end{equation}
where $D$ is the first derivative matrix (Toeplitz) and $\Theta = [\theta(x_0), \theta(x_1), ...,  \theta(x_{N-1})]$
with $\theta(x_i) = 1 - \chi(x_i) + \eta \chi(x_i))$ and $F= [f(x_0), f(x_1), ...,  f(x_{N-1})]$ are vectors in $\R^N$.

The matrix $A = -D \Theta(x) D$ is singular (it has an eigenvalue $0$) and a solution only exists if $F$ is in the image of $A$.
For solving the linear system thus special care has to be taken using either the pseudoinverse, or removing one equation.
This point will be addressed later. 

\medskip

\begin{figure}[htb]
\centering
\includegraphics[width=0.45\textwidth, height= 4.5cm]{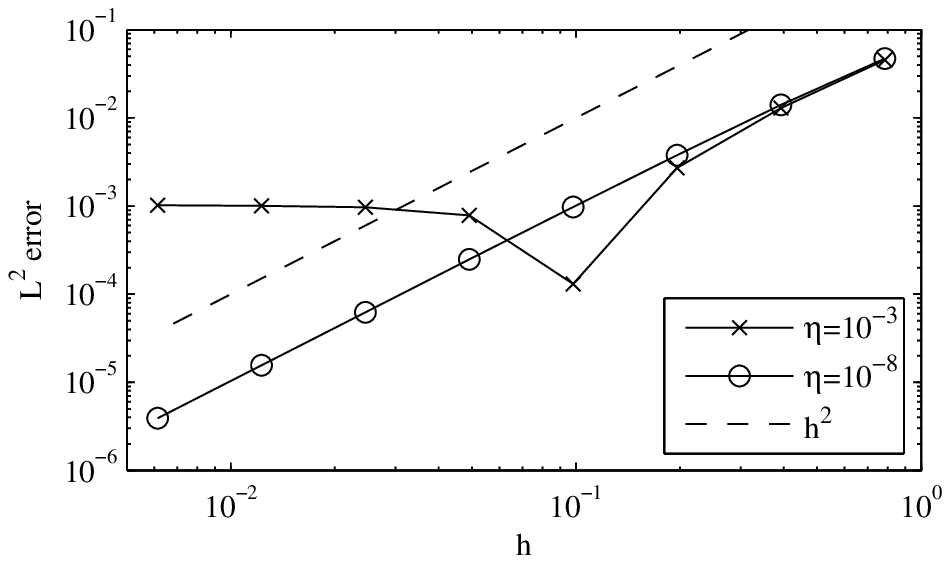}
\includegraphics[width=0.45\textwidth, height= 4.5cm]{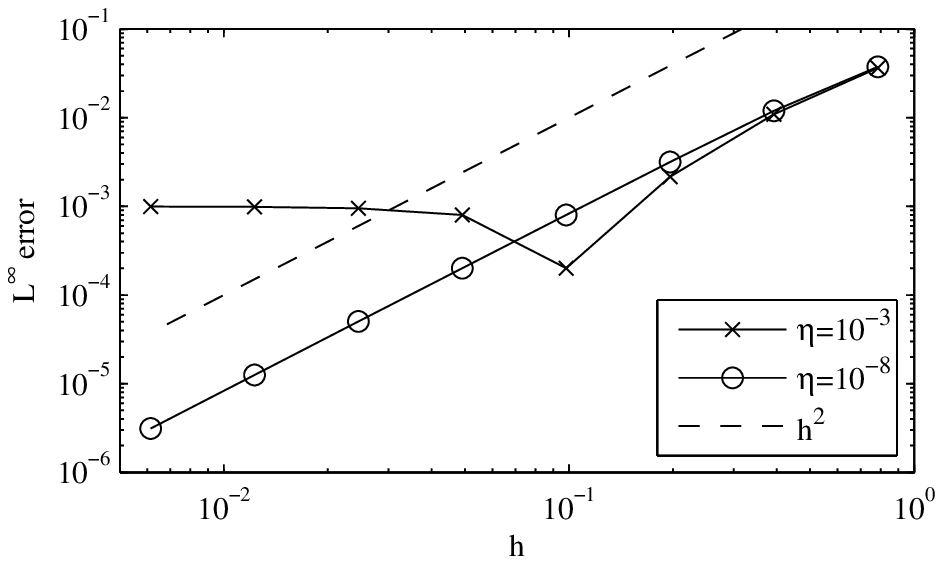}
\caption{Convergence of the second order finite difference scheme for $m=1$.
The $L^2$ (left) and $L^\infty$ (right) errors are calculated only in the fluid
domain $\Omega_f$.}\label{fig:conv_fd2}
\end{figure}

The penalized differential operator can then be approximated to the
second order accuracy with the following finite-difference scheme:
\begin{equation}
 A = -\frac{1}{2} \left( D_F \Theta(x) D_B + D_B \Theta(x) D_F
 \right),
\label{eq:fd2_operator}
\end{equation}
where $D_B$ and $D_F$ are the backward and forward first
derivative matrices,
\begin{equation}
  D_{B} = \frac{1}{h} \left(\begin{array}{ccccc}
    1  &   &   & & -1 \\
    -1 & 1 &   & & \\
    ~  &   &   & \ddots & \\
       &   &   & -1 & 1 \\
  \end{array}\right),
\quad\quad
  D_{F} = \frac{1}{h} \left(\begin{array}{ccccc}
    -1 & 1  &   & & \\
       & -1 & 1 & & \\
    ~  &    &   & \ddots & \\
    1  &    &   & & -1 \\
  \end{array}\right)
\label{eq:df_db_fd2}
\end{equation}
where $h=2\pi/N$.
Note that $dim~ker(A)=1$ reflecting the fact that the (periodic)
solution is defined up to an additive constant. We fix this
constant by imposing the mean value to be zero,
\begin{equation}
 F_1=0, ~~ A_{1,j} = 1,~~j=1,...,N,
\label{eq:fd2_zeromean}
\end{equation}
where $N=dim(A)$. This yields an invertible matrix.
Figure~\ref{fig:conv_fd2} confirms the second-order rate of
convergence, provided that $\eta$ is sufficiently small.

Note that we found that defining the mask function (eq.~\ref{eq:mask1d}) using either the value $0$ or $1$ at the interface,
instead of $1/2$, yields very similar results.

\section{Eigenvalue problem of the penalized Laplace operator}

\subsection{Exact eigenvalue problem}

Now we consider the eigenvalue problem of the Laplace operator with homogeneous Neumann boundary conditions,
\begin{equation}
 - \psi'' \, = \, \lambda \psi \quad \quad x \in (0, \pi)
\end{equation}
with $\psi'(0) = \psi'(\pi) = 0$.
The resulting eigenfunctions are $\psi_n(x) = \cos(n x)$ for $n \in \N$ and the
corresponding eigenvalues are given by $\lambda_n = n^2$. 
Typically, the eigenfunctions are normalized with respect to the $L^2$ norm
and thus the factor $\sqrt{2/\pi}$ has to be included and for $n=0$ we have $\psi_0 = 1/\sqrt \pi$.

\subsection{Penalized eigenvalue problem}

The eigenvalue problem of the penalized Laplace operator with homogeneous Neumann boundary conditions reads,
\begin{eqnarray}
- \phi'' \, &=& \, \lambda \phi \quad \quad \mbox{\rm for} \quad x \in ]0, \pi [ \\
- \eta \phi''  \, &=& \, \lambda \phi \quad \mbox{\rm for} \quad x \in ]\pi, 2 \pi [
\end{eqnarray}
where $\eta >0$ and periodic boundary conditions are imposed at $0$ and $2 \pi$.
Imposing continuity of the solution and of the flux, the problem can be solved exactly and we obtain
the eigenfunctions
\begin{equation}
\phi(x) \, = \,   \left \{
    \begin{array} {ll}
        A_1 \cos (\sqrt{\lambda} x) + B_1 \sin (\sqrt{\lambda} x) \quad \quad \mbox{\rm for} \quad \quad 0 < x < \pi \\
        A_2 \cos (\sqrt{\lambda/\eta} x) + B_2 \sin (\sqrt{\lambda/\eta} x) \quad \mbox{\rm for} \quad \quad \pi < x < 2 \pi \\
    \end{array}
  \right.
\end{equation}
where the coefficients are given by solving the linear system
\begin{eqnarray}
A_1 \cos (\sqrt{\lambda} \pi^-) + B_1 \sin (\sqrt{\lambda} \pi^-) = A_2 \cos (\sqrt{\lambda/\eta} \pi^+) + B_2 \sin (\sqrt{\lambda/\eta} \pi^+)\\ 
-A_1 \sin (\sqrt{\lambda} \pi^-) + B_1 \cos (\sqrt{\lambda} \pi^-) = -A_2 \sqrt{\eta} \sin (\sqrt{\lambda/\eta} \pi^+) + B_2 \sqrt{\eta} \cos (\sqrt{\lambda/\eta} \pi^+)\\
A_1 \cos (\sqrt{\lambda} 0^+) + B_1 \sin (\sqrt{\lambda} 0^+) = A_2 \cos (\sqrt{\lambda/\eta} 2\pi^-) + B_2 \sin (\sqrt{\lambda/\eta} 2\pi^-)\\ 
-A_1 \sin (\sqrt{\lambda} 0^+) + B_1 \cos (\sqrt{\lambda} 0^+) = -A_2 \sqrt{\eta} \sin (\sqrt{\lambda/\eta} 2\pi^-) + B_2 \sqrt{\eta} \cos (\sqrt{\lambda/\eta} 2\pi^-)
\end{eqnarray}
The coefficients $A_1$ and $B_1$ can be eliminated and we obtain a homogeneous linear system
for the coefficients $A_2$ and $B_2$.
\begin{eqnarray}
\left(\begin{array}{cc}
    a  & b\\
    c  & d\\
  \end{array}\right)
\left(\begin{array}{c}
    A_2\\
    B_2\\
  \end{array}\right)
\, = \,  
\left(\begin{array}{c}
    0\\
    0\\
  \end{array}\right)
\end{eqnarray}
with coefficients
\begin{eqnarray}
 a  =  \cos (\sqrt{\lambda/\eta} 2\pi^-) \cos (\sqrt{\lambda} \pi^-) - \sqrt{\eta} \sin (\sqrt{\lambda/\eta} 2\pi^) \sin (\sqrt{\lambda} \pi^-) - \cos (\sqrt{\lambda/\eta} \pi^+) \\
 b  =  \sin (\sqrt{\lambda/\eta} 2\pi^-) \cos (\sqrt{\lambda} \pi^-) + \sqrt{\eta} \cos (\sqrt{\lambda/\eta} 2\pi^) \sin (\sqrt{\lambda} \pi^-) - \sin (\sqrt{\lambda/\eta} \pi^+) \\
 c =  - \cos (\sqrt{\lambda/\eta} 2\pi^-) \sin (\sqrt{\lambda} \pi^-) - \sqrt{\eta} \sin (\sqrt{\lambda/\eta} 2\pi^) \cos (\sqrt{\lambda} \pi^-) + \sqrt{\eta} \sin (\sqrt{\lambda/\eta} \pi^+) \\
 d  = - \sin (\sqrt{\lambda/\eta} 2\pi^-) \sin (\sqrt{\lambda} \pi^-) + \sqrt{\eta} \cos (\sqrt{\lambda/\eta} 2\pi^) \cos (\sqrt{\lambda} \pi^-) - \sqrt{\eta} \cos (\sqrt{\lambda/\eta} \pi^+)
\end{eqnarray}
The eigenvalues $\lambda$ can then be determined by computing the zeros of the determinant of the linear system, {\it i.e.},
solving the nonlinear equation
\begin{equation}
G(\lambda; \eta) \, = \, a d - b c \, = \, 0 
\label{eq:ev_function}
\end{equation}
for a given value of $\eta$.
We did not succeed solving this system symbolically for arbitrary $\eta$, but we can make the following observations:
\begin{itemize}
\item The function $G$ is a periodic function in $\sqrt{\lambda/\eta}$.
\item The value $\lambda=0$ is a solution of eq.(~\ref{eq:ev_function}) and thus an eigenvalue of the penalized operator.
\item The values $\lambda = i^2$ and $\lambda= \eta i^2$ for $i \in \N$ play a special role as different terms in eq.(~\ref{eq:ev_function}) vanish.
\item For the special choice of the penalization parameter $\eta = i^2/j^2$ with $i,j \in \N$, we have explicit solutions and the eigenvalues are $\lambda = i^2$ and $\lambda= \eta i^2$, for $i \in \N$.
\end{itemize}
The above findings motivate the fact that $\lambda = i^2$ and $\lambda= \eta i^2$ are indeed good approximations of the zeros of $G$
for general values of $\eta \in \R^+$.

\subsection{Numerical solution of the penalized eigenvalue problem}

The penalized eigenvalue problem is now solved numerically using second order finite differences.
Thus we discretize,
\begin{equation}
- d_x (\theta(x)) d_x u  \, = \, \lambda u \quad \mbox{\rm for} \quad x \in (0, 2 \pi)
\label{poisson1d_evalue}
\end{equation}
using eq.~(\ref{eq:fd2_operator}) where periodic boundary conditions are imposed at $0$ and $2 \pi$.
The operator $- d_x (\theta(x)) d_x$ is self-adjoint and semi-positive definite, hence all eigenvalues $\lambda$ are real and
positive. 

\bigskip

\begin{figure}[htb]
\centering
\includegraphics{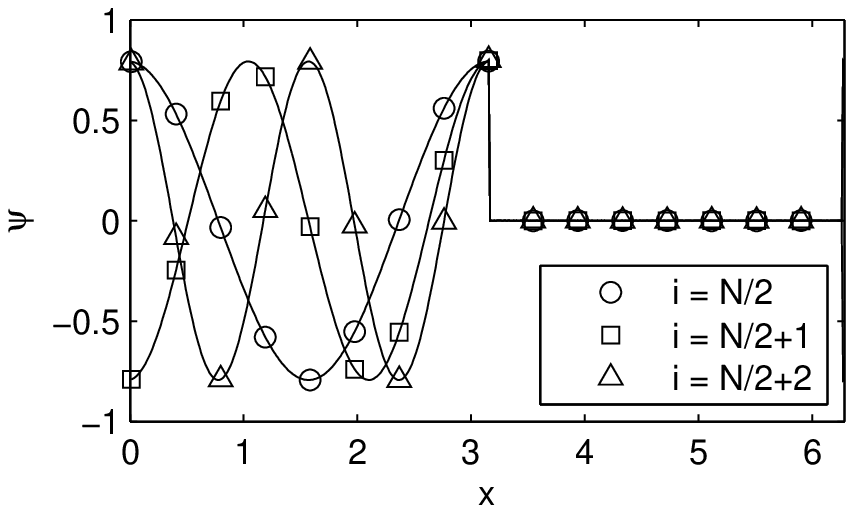}
\caption{Eigenfunctions number $N/2$, $N/2+1$ and $N/2+2$ of the finite-difference penalized Laplace operator. $N=512$, $\eta = 10^{-8}$.}\label{fig:eigenvectors}
\end{figure}

\begin{figure}[htb]
\centering 
\includegraphics[width=0.48\textwidth, height=4.5cm]{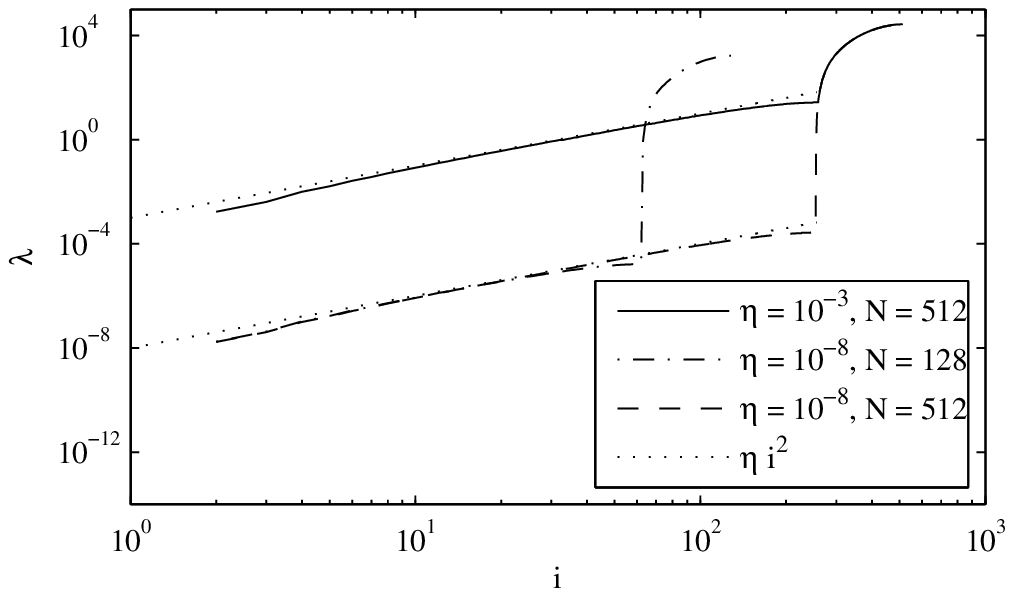}
\includegraphics[width=0.48\textwidth, height=4.5cm]{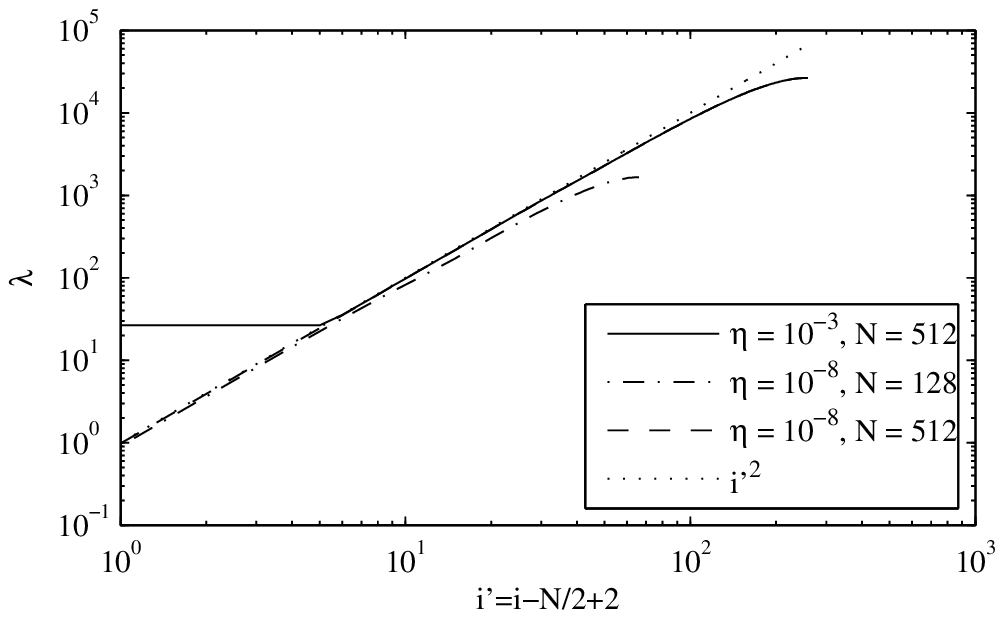}
\caption{Double-logarithmic plot of the eigenvalues $\lambda_i$.
Left: Eigenvalues sorted by their magnitude, in the ascending order. The zero eigenvalue is not shown because of the logarithmic scale.
Right: Eigenvalues in the upper half of the spectrum correspond to the physically relevant ones.}\label{fig:eigenvalues}
\end{figure}

The finite-difference penalized Laplace operator has also a zero eigenvalue, since the solution of the boundary-value problem is only defined up to an additive constant.
One can also identify eigenfunctions of the penalized problem that correspond to the eigenmodes of the original boundary-value problem.
Three of them are displayed in figure~\ref{fig:eigenvectors}.
They correspond to eigenvalues number $N/2$, $N/2+1$ and $N/2+2$.
In the fluid domain (or physical domain, or low-diffusivity domain) they behave like $\cos n x$, and they are close to zero in the other half of the domain.
Similar eigenfunctions exist in the solid (fictitious domain, or large-diffusivity domain), they correspond to the largest eigenvalues.
All non-zero eigenvalues sorted by their magnitude, in the accending order, are shown in figure~\ref{fig:eigenvalues}
for three choices of the model parameters: $N=512$, $\eta=10^{-3}$, $N=128$, $\eta=10^{-8}$ and $N=512$, $\eta=10^{-8}$. 
The spectrum $\lambda_i$ changes from an $\eta i^2$ power law to a concave function approximately at $i=N/2$ (figure~\ref{fig:eigenvalues}, left).
Applying a shift ($i' = i - N/2 +2$) and replotting the upper half of the spectrum for $i \geq  N/2-1$ shows again a power law behavior $\propto i^2$
as illustrated in (figure~\ref{fig:eigenvalues}, right).
For increasing resolution $N$, we can observe that these eigenvalues in the upper half of the spectrum do indeed converge versus the eigenvalues of
the non-penalized Laplace operator given by $i^2$.
The eigenvalues in the lower part of the spectrum depend on the penalization parameter $\eta$ and do converge to zero for $\eta \rightarrow 0$.

\begin{figure}[htb!]
\centering
\includegraphics[width=0.45\textwidth]{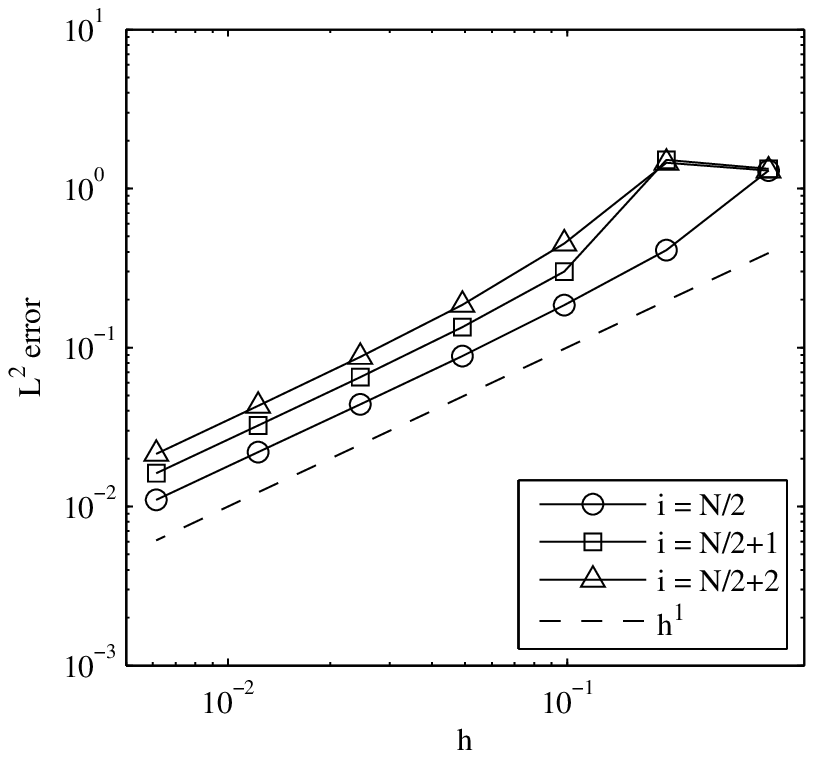}
\quad 
\includegraphics[width=0.45\textwidth]{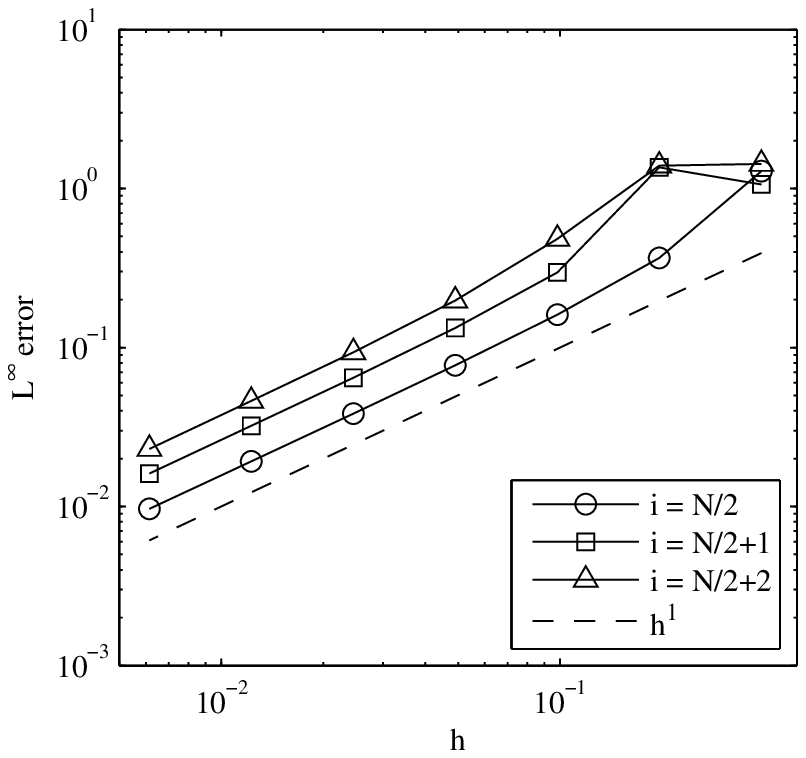}
\caption{$L^2$ (left) and $L^\infty$ (right) distance between the $N/2$-th, $N/2+1$-th and $N/2+2$-th eigenfunctions of the discrete penalized Laplace operator and 2nd, 3rd and 4th eigenfunctions of the continuous Laplace operator with Neumann boundary conditions.}
\label{fig:eig_conv}
\end{figure}

The upper half of the spectrum corresponds to the modes that are only non-trivial in either part of the domain (despite some small oscillations), like in figure~\ref{fig:eigenvalues}.
The lower half of the spectrum corresponds to modes that oscillate with the grid frequency in either subdomain.
Figure~\ref{fig:eig_conv} shows the decay of the distance between the eigenfunctions of the discrete penalized operator (like those in figure~\ref{fig:eigenvectors}) and their exact counterparts, as $h$ decreases. 
In this example, the penalization parameter $\eta=10^{-8}$ is sufficiently small so that the penalization error is smaller than the discretization error within the range of $h$ shown in the figure.
These computations suggest that the discrete eigenfunctions considered here are only a first-order approximation to those of the original boundary-value problem, whereas (we remind that) the solution to the Poisson equation is second-order accurate in $h$.

\section{Application to the penalized Poisson equation in 2d}

Now, we consider a Poisson equation in two space dimensions complemented with homogeneous Neumann boundary conditions,
$$-\nabla^2 u = f$$ with $\partial_n u = 0$.
First, we consider a square domain and then a circular domain. 

\bigskip

The two-dimensional penalized equation in Cartesian coordinates reads
\begin{equation}
- \partial_x (\theta(x,y) \partial_x u(x,y)) - \partial_y (\theta(x,y) \partial_y u(x,y)) = f(x,y).
\label{eq:lap_pen_2d}
\end{equation}
The partial derivatives are approximated using the same second order finite-difference scheme that led to (\ref{eq:fd2_operator}).

\begin{figure}[h!]
\centering
\includegraphics{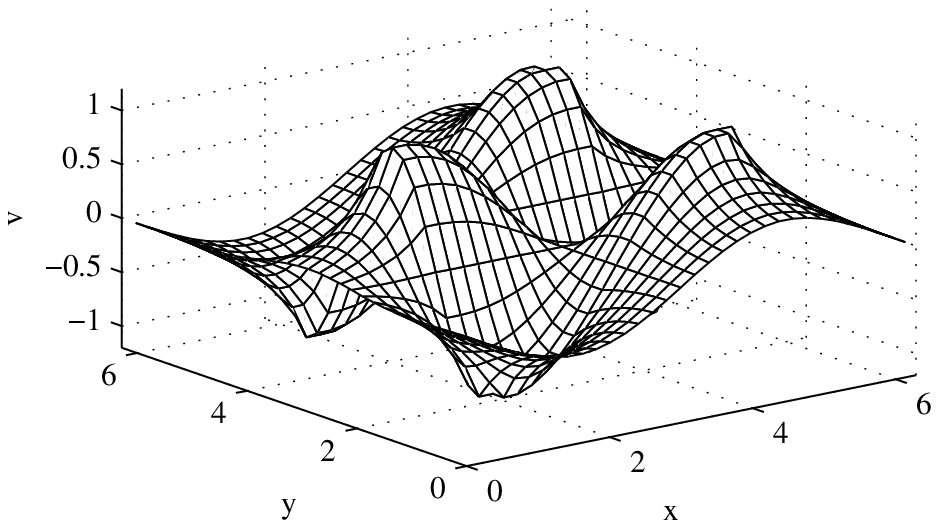}
\caption{Numerical solution of the two-dimensional penalized equation (\ref{eq:lap_pen_2d}) with $\eta=10^{-8}$ and $N=32$ in a rectangular domain.}
\label{fig:example1_2d}
\end{figure}

\begin{figure}[h!]
\centering
\includegraphics{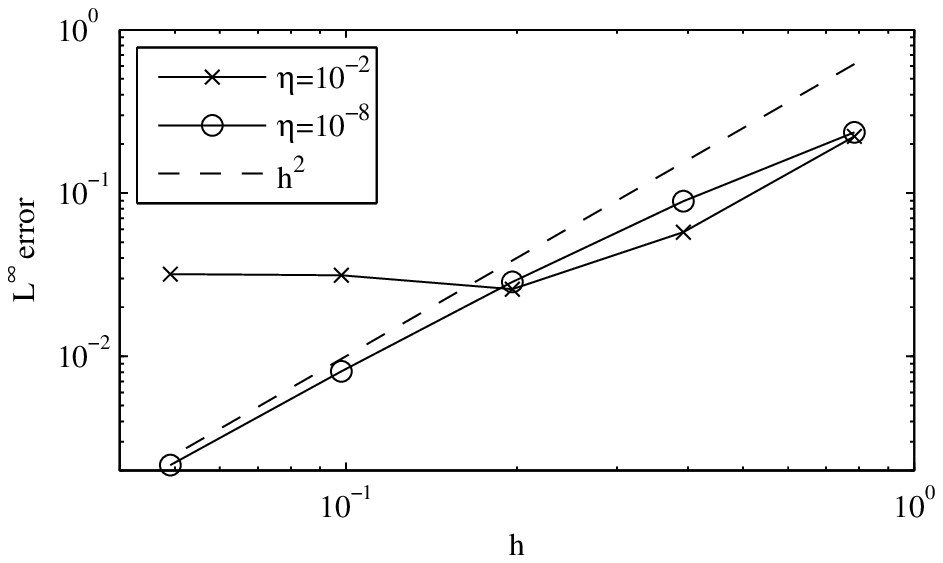}
\caption{$L^\infty$-error decay of the numerical solution of (\ref{eq:lap_pen_2d}) with respect to the exact solution (\ref{eq:exact_2d_ex1}) of the Poisson equation with Neumann boundary conditions in a rectangular domain.
$h=2\pi/N$ is the discretization step size and $\eta$ is the penalization parameter.}\label{fig:conv_2d_ex1}
\end{figure}

Let us first consider an example in which the interface is aligned with the grid.
The computational domain is a periodization of a square $\Omega = [0,2\pi] \times [0,2\pi]$, and the fluid occupies a smaller square sub-domain, $\Omega_f = [\pi/2, 3\pi/2] \times [\pi/2, 3\pi/2]$.
Thus, the mask function is
\begin{equation}
\chi(x,y) = \left\{
\begin{array}{ll}
0 & \textrm{if } x \in ]\frac{\pi}{2},\frac{3\pi}{2}[ ~\textrm{and}~ y \in ]\frac{\pi}{2},\frac{3\pi}{2}[;\\
\frac{1}{2} & \textrm{if } x = \frac{\pi}{2}, y \in ]\frac{\pi}{2},\frac{3\pi}{2}[ ~\textrm{or}~ x = \frac{3\pi}{2}, y \in ]\frac{\pi}{2},\frac{3\pi}{2}[ \\ 
 & ~\textrm{or}~ y = \frac{\pi}{2}, x \in ]\frac{\pi}{2},\frac{3\pi}{2}[ ~\textrm{or}~ y = \frac{3\pi}{2}, x \in ]\frac{\pi}{2},\frac{3\pi}{2}[; \\
\frac{1}{4} & \textrm{if } x = \frac{\pi}{2}, y = \frac{\pi}{2} ~\textrm{or}~ x = \frac{3\pi}{2}, y = \frac{\pi}{2} \\
 & ~\textrm{or}~ x = \frac{\pi}{2}, y = \frac{3\pi}{2} ~\textrm{or}~ x = \frac{3\pi}{2}, y = \frac{3\pi}{2}; \\
1 & \textrm{otherwise}
\end{array}\right.
\label{eq:mask_2d_ex1}
\end{equation}
Let the right-hand side of the penalized Poisson equation (\ref{eq:lap_pen_2d}) be
\begin{equation}
f(x,y) = \left\{
\begin{array}{ll}
5 \sin x \cos 2y & \textrm{if } x \in ]\frac{\pi}{2},\frac{3\pi}{2}[ ~\textrm{and}~ y \in ]\frac{\pi}{2},\frac{3\pi}{2}[;\\
\frac{5}{2} \cos 2y & \textrm{if } x = \frac{\pi}{2}, y \in ]\frac{\pi}{2},\frac{3\pi}{2}[; \\
-\frac{5}{2} \cos 2y & \textrm{if } x = \frac{3\pi}{2}, y \in ]\frac{\pi}{2},\frac{3\pi}{2}[; \\
-\frac{5}{2} \sin x & \textrm{if } y = \frac{\pi}{2}, x \in ]\frac{\pi}{2},\frac{3\pi}{2}[ ~\textrm{or}~ y = \frac{3\pi}{2}, x \in ]\frac{\pi}{2},\frac{3\pi}{2}[; \\
-\frac{5}{4} & \textrm{if } x = \frac{\pi}{2}, y = \frac{\pi}{2} ~\textrm{or}~ x = \frac{\pi}{2}, y = \frac{3\pi}{2}; \\
\frac{5}{4} & \textrm{if } x = \frac{3\pi}{2}, y = \frac{\pi}{2} ~\textrm{or}~ x = \frac{3\pi}{2}, y = \frac{3\pi}{2}; \\
0 & \textrm{otherwise}
\end{array}\right.
\label{eq:rhs_2d_ex1}
\end{equation}
Note that the zero mean value of the numerical solution in the fluid domain is imposed by replacing the first equation in the linear system by
\begin{equation}
\sum_{i,j=\overline{1,N}}\left[1-\chi(x_{ij})\right]u_{ij} = 0.
\label{eq:mean_fluid}
\end{equation}
In the fluid domain $\Omega_f$, the solution to (\ref{eq:lap_pen_2d}) converges to 
\begin{equation}
w(x,y) = \sin x \cos 2y \textrm{,    where  } x \in ]\frac{\pi}{2},\frac{3\pi}{2}[ ~\textrm{and}~ y \in ]\frac{\pi}{2},\frac{3\pi}{2}[\,,
\label{eq:exact_2d_ex1}
\end{equation}
as $\eta \to 0$.
Figure~\ref{fig:example1_2d} displays a numerical solution to this problem with $N=32$ discretization grid points in each direction and with $\eta=10^{-8}$.
Inside the fluid domain, the solution is close to (\ref{eq:exact_2d_ex1}).
Outside, it is close to a harmonic function (up to numerical errors).
Figure~\ref{fig:conv_2d_ex1} shows the decay of the $L^\infty$ error of the finite-difference solution 
with respect to the exact solution (\ref{eq:exact_2d_ex1}) in the fluid domain (including the points on the boundary).
Two values of $\eta$ are considered. For $\eta=10^{-2}$, the error saturates at $h < 0.2$, where $h=2\pi/N$. 
For $\eta=10^{-8}$, the decay approaches the theoretical $-2$ slope for small $h$ and the saturation 
is not observed within this range of $h$, implying that the penalization error is much smaller than the discretization error.

\bigskip

Let us consider a circular fluid domain,
with the mask function
\begin{equation}
\chi(x,y) = \left\{
\begin{array}{ll}
0 & \textrm{if } r < \pi;\\
\frac{1}{2} & \textrm{if } r = \pi; \\ 
1 & \textrm{otherwise,}
\end{array}\right.
\label{eq:mask_2d_ex2}
\end{equation}
where $r=\sqrt{(x-\pi)^2+(y-\pi)^2}$. 
The right-hand side is
\begin{equation}
f(x,y) = \left\{
\begin{array}{ll}
4\cos2r + \frac{2\sin2r}{r} & \textrm{if } r < \pi;\\
-\frac{1}{2} & \textrm{if } r = \pi; \\ 
0 & \textrm{otherwise.}
\end{array}\right.
\label{eq:rhs_2d_ex2}
\end{equation}
The exact solution to the Poisson equation with homogeneous Neumann boundary conditions in this case is
\begin{equation}
w = \cos 2r + \frac{4}{\pi^2} \textrm{,    where  } r < \pi\,,
\label{eq:exact_2d_ex2}
\end{equation}
inside the fluid domain embedded in a square computational domain $\Omega=[0,2\pi]\times[0,2\pi]$.

\begin{figure}[htb]
\centering
\includegraphics[scale=0.15]{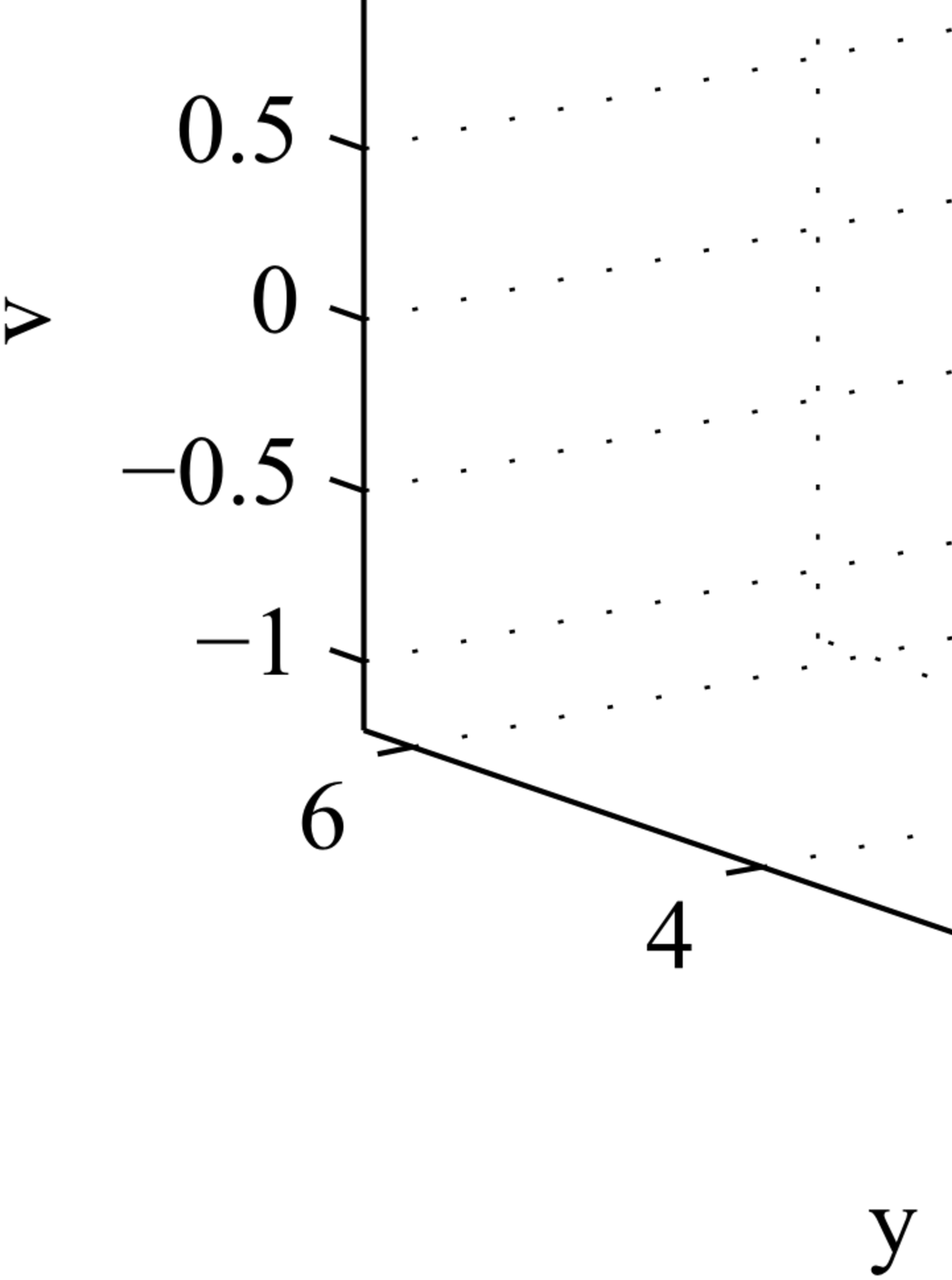}
\quad 
\includegraphics[scale=0.15]{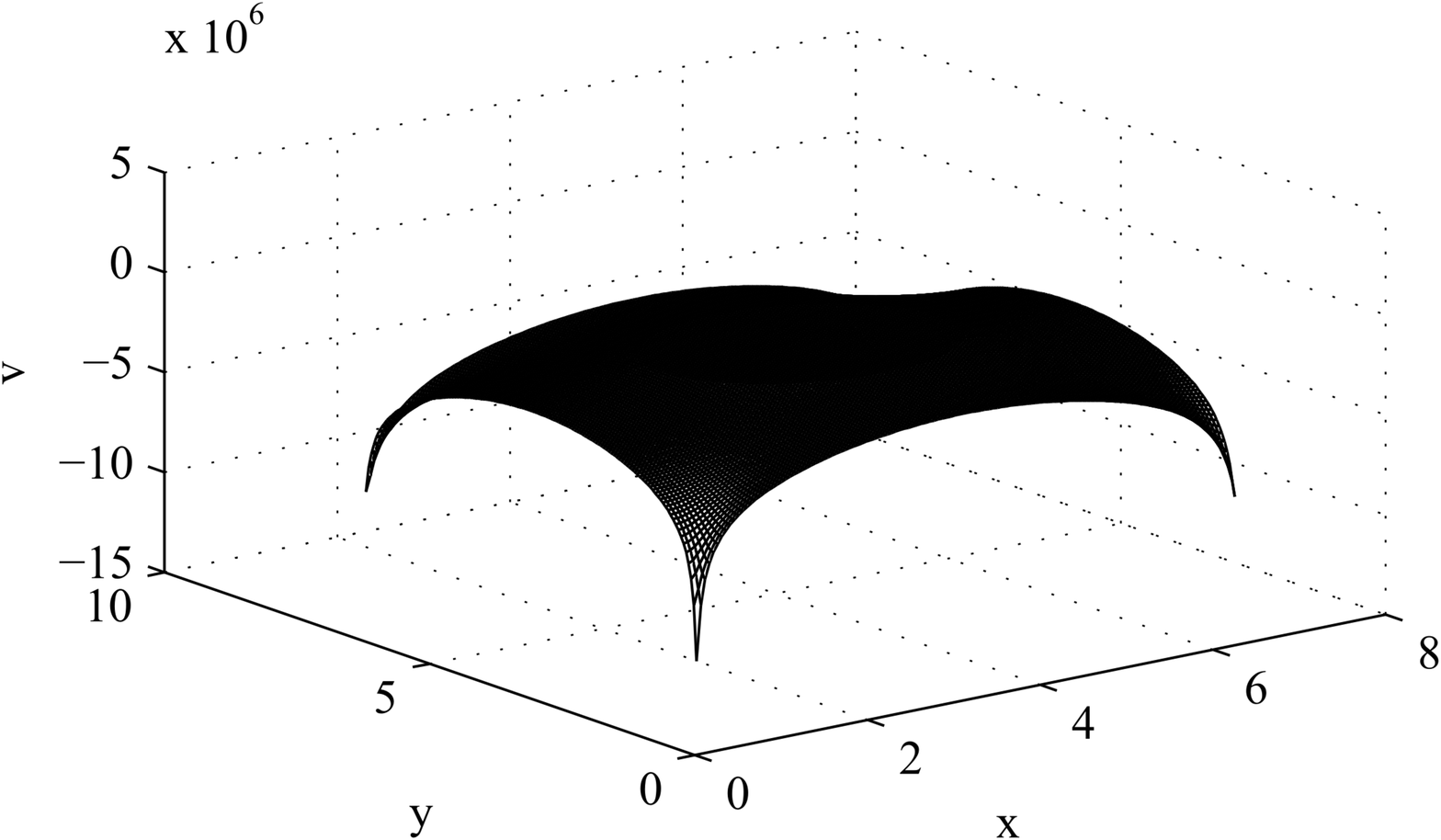}
\caption{Numerical solution of the two-dimensional penalized equation (\ref{eq:lap_pen_2d}), (\ref{eq:mask_2d_ex2}), (\ref{eq:rhs_2d_ex2}) 
in a circular domain with $\eta=10^{-8}$ and $N=127$, first linear equation replaced with the zero-mean condition. Top: Zoom of solution in the fluid domain.  Bottom: Total domain
illustrating the singular behavior in the solid domain.}
\label{fig:example2_2d_1}
\end{figure}

\begin{figure}[htb]
\centering
\includegraphics[scale=0.15]{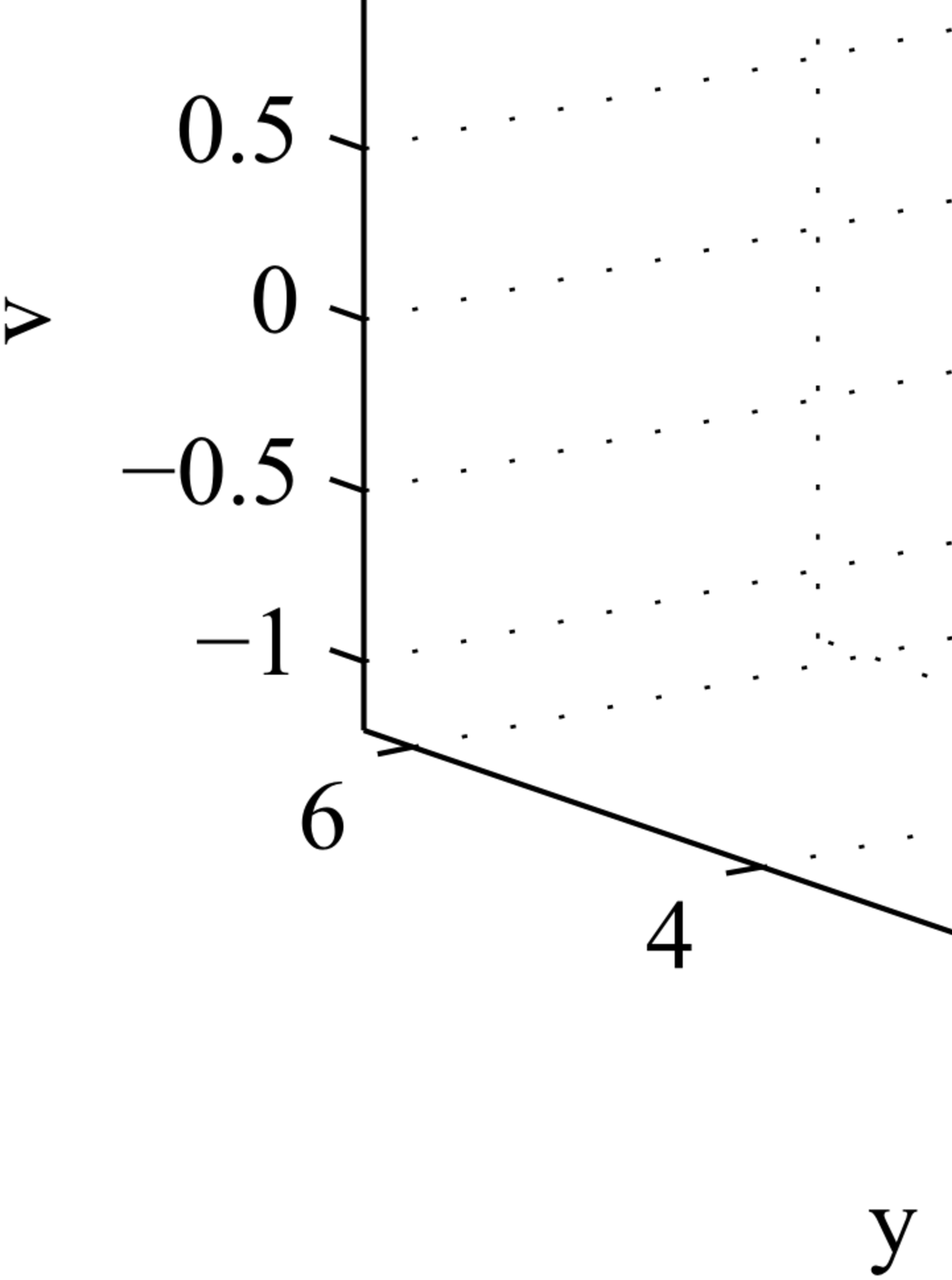}
\quad 
\includegraphics[scale=0.15]{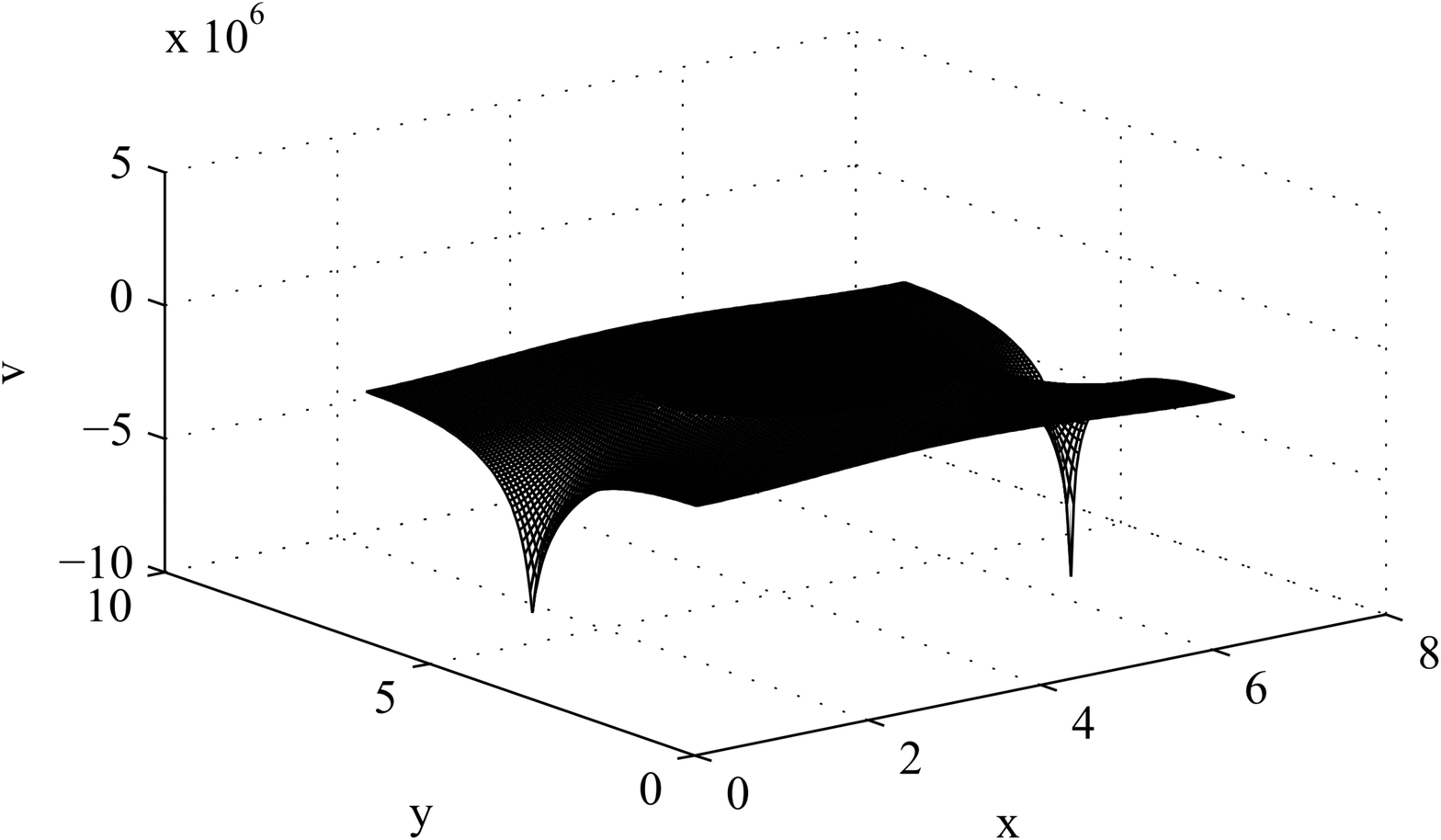}
\caption{Numerical solution of the two-dimensional penalized equation (\ref{eq:lap_pen_2d}), (\ref{eq:mask_2d_ex2}), (\ref{eq:rhs_2d_ex2}) 
in a circular domain with $\eta=10^{-8}$ and $N=127$, $N^2/2$-th linear equation replaced with the zero-mean condition. Top: Zoom of solution in the fluid domain.  Bottom: Total domain
illustrating the singular behavior in the solid domain.}
\label{fig:example2_2d_2}
\end{figure}

\begin{figure}[htb]
\centering
\includegraphics[scale=0.15]{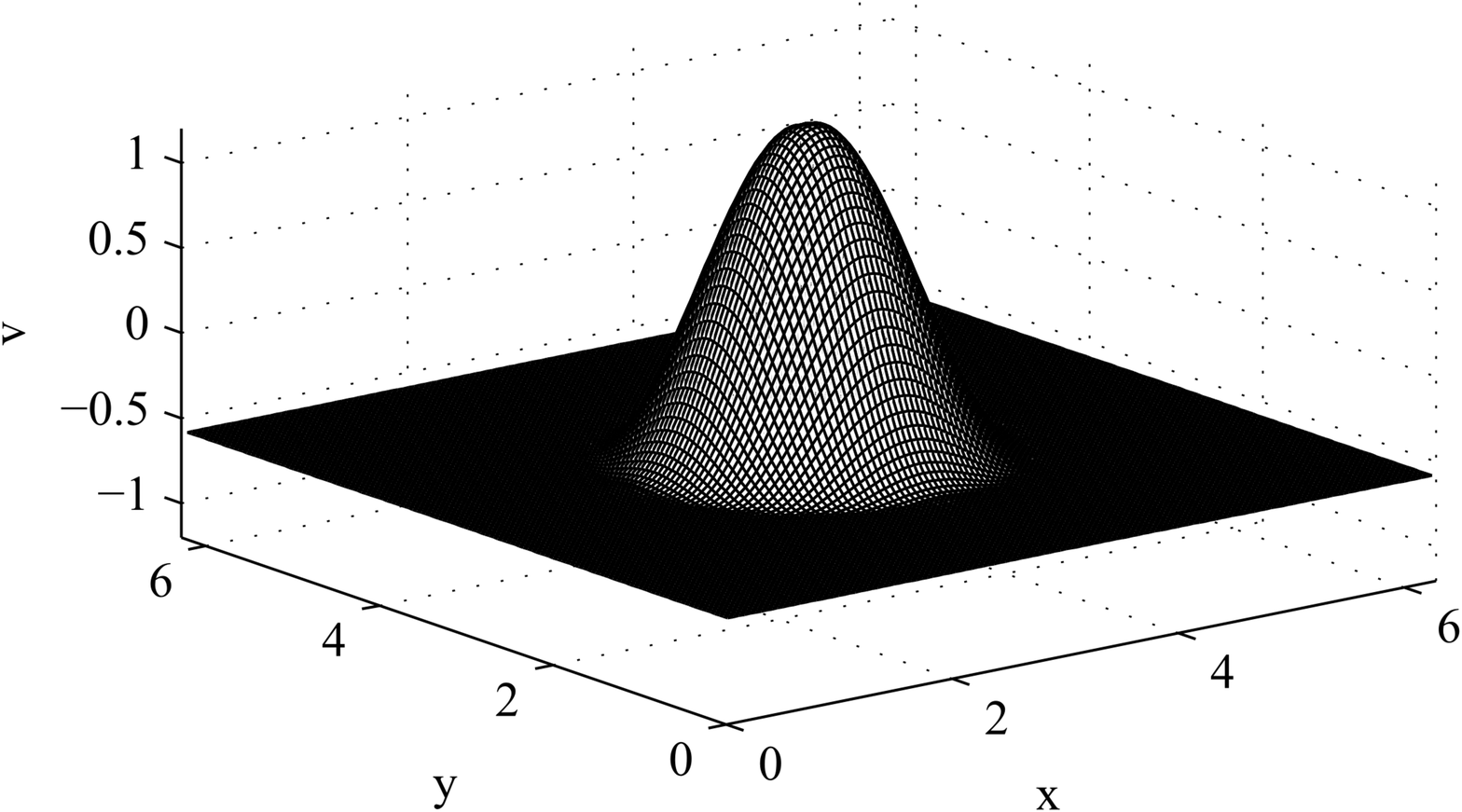}
\quad 
\includegraphics[scale=0.15]{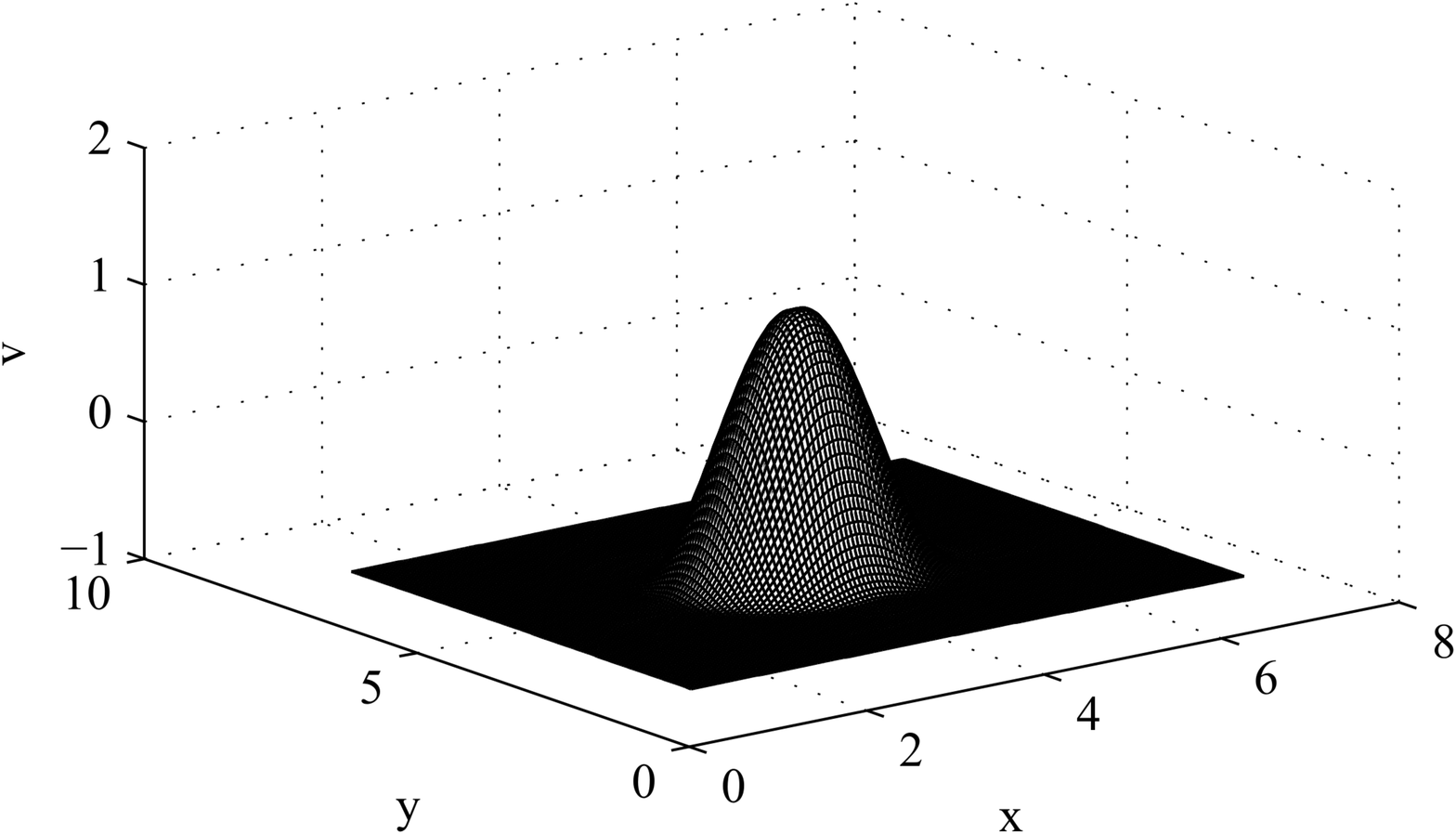}
\caption{Numerical solution of the two-dimensional penalized equation (\ref{eq:lap_pen_2d}), (\ref{eq:mask_2d_ex2}), (\ref{eq:rhs_2d_ex2}) 
in a circular domain with $\eta=10^{-8}$ and $N=127$, $N^2/2+N/2$-th linear equation replaced with the zero-mean condition. Top: Zoom of solution in the fluid domain.  Bottom: Total domain illustrating the smooth behavior in the solid domain.}\label{fig:example2_2d_3}
\end{figure}

\begin{figure}[h!]
\centering
\includegraphics{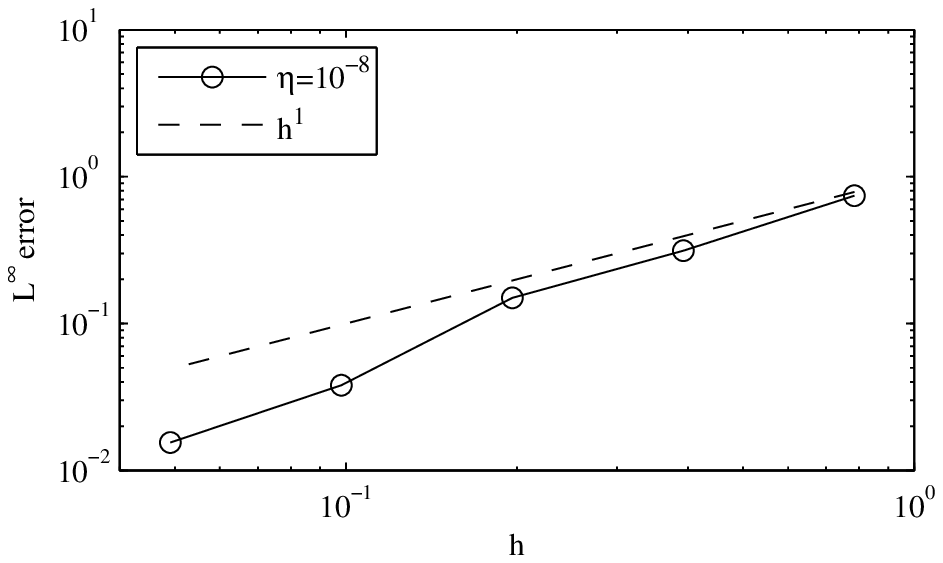}
\quad
\includegraphics{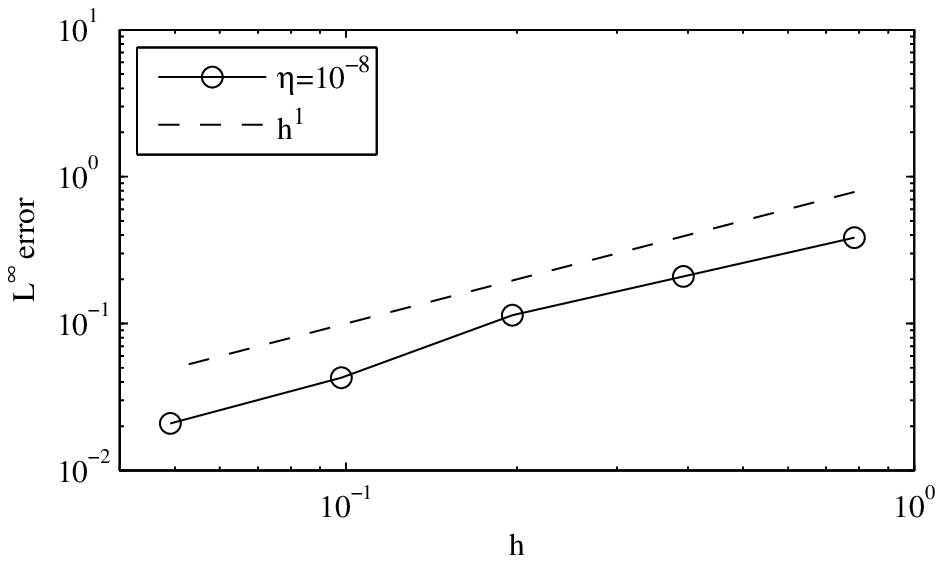}
\quad
\includegraphics{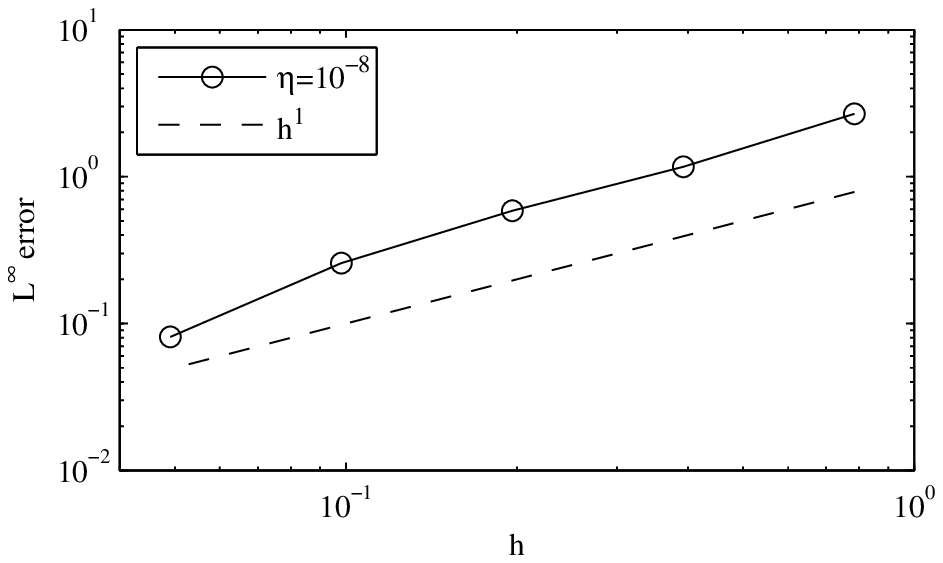}
\caption{Convergence plots of the two-dimensional penalized equation in a circular domain for the three above cases (see figures~\ref{fig:example2_2d_1}, \ref{fig:example2_2d_2} and \ref{fig:example2_2d_3}, respectively).}\label{fig:conv_2d_ex2}
\end{figure}

We observed that the numerical solution is sensitive to the choice of the linear equation which is replaced with the zero-mean condition.
The operator matrix has many small eigenvalues if $\eta$ is small.
Another possibility would be to add the zero-mean condition without removing any of the equations and solve an overdetermined system in the least-square sense (results not shown here). Note that in this case we observed a smooth behavior in the solid domain.
Figure~\ref{fig:example2_2d_1} shows the solution with the first equation replaced and $\eta=10^{-8}$, $N=127$.
Figure~\ref{fig:example2_2d_2} displays the same solution with the $N^2/2$-th equation replaced, and figure~\ref{fig:example2_2d_3} with the $N^2/2+N/2$-th equation replaced.
The solution in the fluid is slightly different in the three cases (and seems to be convergent with $\eta$ and $h=2\pi/N$).
In the solid domain, a parasite harmonic solution appears, which has a singularity at the point that corresponds to the removed equation.
The convergence of the two-dimensional penalized equation for the three above cases is summarized in Figure~\ref{fig:conv_2d_ex2}
and shows first order convergence in all cases.
The second order convergence observed in the one-dimensional case (subsec. \ref{subsec_defds}) and for the two-dimensional case in the rectangular domain is thus reduced to first order. The reason is that the Cartesian grid introduces a staircase effect and 
the approximation of the circular mask function reduces to first order. 
Techniques to obtain higher order for complex geometries (based on interpolation) have been proposed in~\cite{SVCA08}.

%

\section{Conclusions}
The volume penalization method to impose homogeneous Neumann boundary conditions
has been analyzed by considering the Poisson equation. 
In one space dimension, the penalized Poisson equation
has been solved analytically for a particular right hand side and the penalization error has been determined showing $O(\eta)$ convergence of the solution towards the solution of the exact problem.
We also found that no penalization boundary layer is present.
This observation is in contrast to what was found for
the time-dependent heat equation with Neumann conditions \cite{KKAS12} and also for the Poisson equation with Dirichlet boundary conditions \cite{NKS12}.
In both cases, there is a penalization boundary layer which becomes thinner for decreasing penalization parameter $\eta$ and its thickness scales like $O(\sqrt{\eta})$.
This implies that only an $O(\sqrt{\eta})$ convergence can be proven \cite{ABF99,CaFa03,KKAS12}. 
Nevertheless for the penalized Laplace operator with Neumann conditions, the corresponding matrix becomes ill--conditioned
and the condition number behaves like $O(1/\eta)$. 
Thus, special care has to be taken for the numerical solution, as in addition to the singularity of the matrix (the presence of an eigenvalue $0$), 
the linear system becomes stiff.

The performed numerical simulations using second order finite differences yield second order convergence of the solution towards the
solution of the Poisson equation, given that the penalization parameter is sufficiently small.
Due to the regularity of the exact solution of the penalized equation and the $O(\eta)$ behavior of the penalization error,
we anticipate that for higher order numerical methods we will also find second order convergence.

The eigenvalue problem of the penalized Laplace operator with Neumann boundary conditions was also studied in some detail.
We found that the spectrum of the penalized operator exhibits two distinct behaviors. 
The upper part of the spectrum corresponding to the large
eigenvalues converges for increasing resolution $N$ to the spectrum of the exact operator ($\propto i^2$).
For the lower part, corresponding to the small eigenvalues, the spectrum exhibits the same power law scaling but
the values are multiplied with $\eta$ and thus converge to zero for $\eta \rightarrow 0$.
The eigenfunctions in the upper half of the spectrum are non-trivial in either part of the domain.
The lower half of the spectrum corresponds to modes that oscillate with the grid frequency in either subdomain.
The eigenfunctions corresponding to the upper half of the spectrum of the discrete penalized operator converge to 
their exact counterparts and we found first order convergence using second order finite differences. 

In two space dimensions, we performed numerical simulations for a rectangular geometry for which
the grid is aligned with the boundary.
In this case we obtained again second order convergence of the numerical solution.
For the circular geometry, for which the boundary is not aligned with the Cartesian grid,
only first order convergence can be observed which is due to the geometrical error.

An interesting perspective is the extension of the volume penalization to higher order penalization, 
also called active penalization, using, {\it e.g.}, smooth extensions of the solution, based for instance on Hermite interpolation, as proposed in \cite{MLBS12}.
First promising results using active penalization for Navier--Stokes have been presented in \cite{ShNa13}. 
%
%
An extension to impose inhomogeneous Neumann conditions has been proposed in \cite{MLBS12}
for Fourier spectral methods. The underlying idea is to use volume penalization to impose
Dirichlet boundary conditions for the derivative and then integrating the equation, which
can be easilty done in spectral space. 

\section*{Acknowledgements}
RNVY is grateful to the Humboldt Foundation for its support through a post-doctoral grant.
KS thanks the organizers of WONAPDE 2013 for their kind invitation to Concepcion, Chile.


\end{document}